\documentclass[12pt]{article}
\usepackage{amssymb}

\title{\sc Differential Calculus for Dirichlet forms:
The measure-valued gradient
 preserved by image
\footnote{This work was presented at the meeting on Stochastic Analysis
and Potential Theory, St Priest de Gimel, 1-6 sept 2002}}
 \author{Nicolas BOULEAU\\ Ecole des Ponts et Chauss\'ees, Paris}

\begin{document}
\maketitle

\begin{abstract}
In order to develop a differential calculus for error propagation (cf [3]) we study local Dirichlet forms on  probability spaces
 with carr\'e du champ $\Gamma$ -- i.e. error structures -- and we are looking for an object related to $\Gamma$ which
is linear and with a good behaviour by images. 
For this we introduce a new notion called the measure valued gradient which is a randomized square root of 
$\Gamma$. The exposition begins with  inspecting some natural notions candidate to 
solve the problem before proposing the measure-valued gradient and proving its satisfactory properties.
\end{abstract}

\section{Preamble}
\label{}
Our main purpose being to study images, in order to avoid unessential difficulties, we restrict us to Dirichlet forms defined on probability spaces.
On a probability space $(W,\;{\mathcal W},\;m)$ let us consider a local Dirichlet form $(\mathbb{D}, \mathcal E)$ with carr\'e du champ operator $\Gamma$.
This is equivalent  (cf.[4], [6]) to the data of 

(1) a dense sub-vector space $\mathbb{D}$ of  $L^2(W,\;{\mathcal W},\;m)$,

(2) a symmetric positive bilinear operator $\Gamma$ from $\mathbb{D}\times\mathbb{D}$ into $L^1(m)$
satisfying the following functional calculus :
\noindent if $u\in\mathbb{D}^m,\;v\in\mathbb{D}^n,\;$ and $ F,G$ are Lipschitz and $C^1$ from $\mathbb{R}^m$ [resp. $\mathbb{R}^n$] into $\mathbb{R}$
 then $F(u)\in\mathbb{D}$ and $G(v)\in\mathbb{D}$ and
$$\Gamma[F(u),G(v)]=\sum_{i,j}F_i^\prime(u)G_j^\prime(v)\Gamma[u_i,v_j]\quad\quad m{\mbox{-a.e.}}$$

(3) and such that the form ${\mathcal E}$ given by $${\mathcal E}[u,v]=\frac{1}{2}\int\Gamma[u,v]\,dm$$ 
 is  {\it closed }  in $L^2(W,{\mathcal W}, m)$, i.e. $\mathbb{D}$ is complete with the norm
$$\|\;.\;\|_\mathbb{D} =\left(\|\;.\;\|_{L^2(m)}^2+{\mathcal E}[\;.\;]\right)^{1/2}.$$

We write always $\Gamma[u]$ for
$\Gamma[u,u]$. A term
$S=(W,\;{\mathcal W},\;m,\;\mathbb{D},\;\Gamma)$
satisfying properties (1) (2) (3) is called an error structure.  
The notion of error structure  is stable 
by products, even infinite products,  and this feature gives easily error structures on spaces of stochastic processes
like the Wiener space (cf. books [4], [6] and [3], and examples of application [1], [2]).

Such a term $(W,\;{\mathcal W},\;m,\;\mathbb{D},\;\Gamma)$ is also easily transported by images :
If $X\in \mathbb{D}^d$,  let us consider the space ${\mathcal C}^1\cap Lip(\mathbb{R}^d,\mathbb{R})$ of functions $u$ of class ${\mathcal C}^1$ and Lipschitz from $\mathbb{R}^d$ into $\mathbb{R}$
(which are such that $u\circ X\in \mathbb{D}$), then the term $S_X=(\mathbb{R}^d,\; {\mathcal B}(\mathbb{R}^d),\; X_\ast m,\; \mathbb{D}_X,\; \Gamma_X)$ 
where $X_* m$ is the image measure of $m$ by $X$, 
$\mathbb{D}_X$ is the closure of ${\mathcal C}^1\cap Lip(\mathbb{R}^d,\mathbb{R})$ for the norm $\| u\|_{\mathbb{D}_X}=\|u\circ X\|_{\mathbb{D}}$ and 
$$\Gamma_X[u](x)=\mathbb{E}[\Gamma[u\circ X]\;|\;X\!=\!x]$$
satisfies still properties (1) (2) (3), i.e. is still an error structure.

The question we attempt to answer here, is to find an object related to $\Gamma$ which be linear and preserved by image.
Let us first look at some existing objects in the literature.

{\it The generator} 
 $(A,\;{\mathcal D}A)$ of the strongly continuous semigroup 
on $L^2(W,\;{\mathcal W},\;m)$ associated
with the error structure is a linear operator and is transformed by image into the generator of the image structure 
by a relation similar to that defining the image of $\Gamma$ :
$$\forall f\;:\; f\circ X\in{\mathcal D}A\quad\quad A_X[f](x)=\mathbb{E}[A[f\circ X]\;|\;X\!=\!x]$$
(cf. [4] chapter V prop. 1.1.7 and 1.1.8)
 but the calculations with $A$ involve non linear operations because of the presence of $\Gamma$:
if $f\in{\mathcal C}^2_{bb}$ and if $\Gamma[X_i,X_j]\in L^2(m)$ we have
$$ A[f\circ X]=\sum_i A[X_i]\frac{\partial f}{\partial x_i}(X)+\frac{1}{2}\sum_{i,j}\Gamma[X_i,X_j]\frac{
\partial^2 f}{\partial x_i\partial x_j}(X).$$

{\it A Dirichlet-gradient} or shortly D-gradient for the error structure \break $(W,\;{\mathcal W},\;m,\;\mathbb{D},\;\Gamma)$ is 
defined with an auxiliary separable Hilbert space $H$ as a linear map $D$
 from $\mathbb{D}$ into $L^2(m,H)$ s.t.

$(i)\quad \forall U\in\mathbb{D}\quad\quad \|D[U]\|^2_H=\Gamma[U].$

Such an operator satisfies necessarily

$(ii)\quad \forall F\in Lip,\;\forall U\in \mathbb{D}\quad D[F\circ U]=F^\prime\circ U.D[U]$

$(iii) \quad \forall F\in {\mathcal C}^1\cap Lip(\mathbb{R}^d),\;\forall U\in \mathbb{D}^d\quad D[F\circ U]=
\sum_i F^\prime_i\circ U.D[U_i]$.

\noindent Any error structure admits a D-gradient as soon as $\mathbb{D}$ is separable (Mokobodzki, cf. [4] exercise 5.9).

Let us emphasize nevertheless that a D-gradient does not give by image a D-gradient for the image structure :

Let $S_X=(\mathbb{R}^d,\; {\mathcal B}(\mathbb{R}^d),\; X_\ast m,\; \mathbb{D}_X,\; \Gamma_X)$ be the image of $S$ by $X\in\mathbb{D}^d$, 
then the formula 
$\mathbb{E}[D[F\circ X]\;|\;X\!=\!x]$ does not define a D-gradient for $X_\ast S$ because
$$
(*)\left\{\begin{array}{l}
\quad<\mathbb{E}[D[F\circ X]\;|\;X\!=\!x]\;,\;\mathbb{E}[D[F\circ X]\;|\;X\!=\!x]>_{_H}\\
\quad=\!\!\!\!\!/\quad\mathbb{E}[<D[F\circ X],D[F\circ X]>_{_H}\;|\;X\!=\!x]\quad\quad(=\Gamma_X[F](x)).
\end{array}\right.
$$

A D-gradient is not a canonical notion, there is  latitude in its definition. The space $H$ in particular may be chosen in different ways depending
on what is the most convenient and simple in the examples.

{\it The Feyel-la-Pradelle derivative} is a particular case of D-gradient which is canonical in
 the case the measure $m$ is Gaussian and $W$ a vector space (cf [5]). It can be generalized to non-Gaussian cases by taking for $H$ the space
$L^2(\hat{W},\hat{\mathcal W}, \hat{m})$ where $(\hat{W},\hat{\mathcal W}, \hat{m})$ is a copy of $(W,{\mathcal W}, m)$ (cf [5] and [4] chapter
III\S 2 and chapter VII\S 1). With respect to our question, it has the same weakness as any D-gradient of being not preserved by images because
of the inequality (*).

\section{The measure valued gradient}
The object we shall define possesses, like the D-gradient, some latitude in its definition and depends on an auxiliary Hilbert space.

We suppose $S=(W,\;{\mathcal W},\;m,\;\mathbb{D},\;\Gamma)$ admits a D-gradient defined with the separable Hilbert space $H$.
For our purpose, let us recall the notion of 
{\it   Gaussian white noise measure}
 based  on a measurable space $(E,{\mathcal F})$ with associated positive measure $\mu$.

\noindent{\bf Definition 1.} {\it  Let $\mu$ be a bounded positive measure on the measurable space $(E,{\mathcal F})$. 
 A Gaussian white noise measure $\nu$ on $(E,{\mathcal F})$ defined on the probability space $(\Omega,{\mathcal A}, \mathbb{P})$
 with  associated positive measure $\mu$ is a map from ${\mathcal F}$ into
$L^2(\Omega,{\mathcal A}, \mathbb{P})$ s. t.

{\rm(j)}\quad$\forall A\in{\mathcal F},\; \nu(A)$ is a centered Gaussian r.v. with variance $\mu(A)$,

{\rm(jj)}\quad $A\mapsto \nu(A)$ is $\sigma$-additive in $L^2(\Omega,{\mathcal A}, \mathbb{P})$,

{\rm(jjj)} \quad if $A_1,\ldots, A_k\in{\mathcal F}$ are pairwise disjoint the r.v. $\nu(A_1),\ldots,\nu(A_k)$ are independent.} 

Then $\nu$ extends uniquely from ${\mathcal F}$ to $L^2(E,{\mathcal F},\mu)$ and we write $\nu(f)$ for $f\in L^2(E,{\mathcal F},\mu)$. Given
 a measured space $(E,{\mathcal F}, \mu)$ such that $L^2(E,{\mathcal F},\mu)$ is separable such a white noise measure may be constructed
as a classical Wiener integral in the following way : let us take $(\Omega,{\mathcal A}, \mathbb{P})=(\mathbb{R},{\mathcal B}(\mathbb{R}), {\mathcal N}(0,1))^\mathbb{N}$
so that the coordinate mappings $(g_n)$ are i.i.d. Gaussian reduced r.v. Then for $f\in L^2(E,{\mathcal F}, \mu)$, we can put
$$\nu(f)=\sum_n (f,\xi_n)_{L^2(\mu)}\,g_n$$ where $(\xi_n)$ is an orthonormal basis of $L^2((E,{\mathcal F}, \mu)$. If $L^2(E,{\mathcal F},\mu)$ is no more 
supposed to be separable, such a white noise measure may be constructed as Gaussian process indexed by $\mathcal F$ by Kolmogorov theorem.

The positive measure $\mu$ associated with the white noise measure $\nu$ will be often denoted by the symbolic notation $\mathbb{E}_{\mathbb{P}}[(d\nu)^2]$.

Similarly, given on $(E,{\mathcal F})$ a symmetric matrix of bounded measures 
$
\left(
\begin{array}{rl}
\mu_{11} &\mu_{12}\\
\mu_{12} &\mu_{22}
\end{array}
\right)
$
such that 
\noindent$
\left(
\begin{array}{rl}
\mu_{11}(A) &\mu_{12}(A)\\
\mu_{12}(A) &\mu_{22}(A)
\end{array}
\right)
$
 be positive $\forall A\in{\mathcal F}$, we can define a bivariate white noise measure which to each $A\in{\mathcal F}$ associates a pair of Gaussian
variables $(\nu_1(A),\nu_2(A))$ satisfying properties analogous to (j), (jj), (jjj). Such a bivariate white noise may be transformed
in different ways :

a) It may be multiplied by a function $\varphi\in L^2(E,{\mathcal F},\mu_{11}+\mu_{22})$ :
$$
(\varphi\nu)(A) =
\left(
\begin{array}{c}
\int_A \varphi\,d\nu_1\\
\int_A \varphi\,d\nu_2
\end{array}
\right)
$$
is a bivariate white noise measure with associated matrix 
$\varphi^2
\left(
\begin{array}{rl}
\mu_{11} &\mu_{12}\\
\mu_{12} &\mu_{22}
\end{array}
\right).
$

b) For $x=(x_1,x_2)\in\mathbb{R}^2$ we can define the scalar white noise measure $(x,\nu)=x_1\nu_1+x_2\nu_2$ whose associated measure
is 
$
x^t\left(
\begin{array}{rl}
\mu_{11} &\mu_{12}\\
\mu_{12} &\mu_{22}
\end{array}
\right)x.
$

c) As a mixing of a) and b) let $\psi=(\psi_1,\psi_2)$ be in $L^2(E,{\mathcal F},\mu_{11}+\mu_{22}; \mathbb{R}^2)$. We can define the scalar
white noise measure $(\psi,\nu)$ by 
$$
(\psi,\nu)(A)=\int_A\psi_1d\nu_1+\int_A\psi_2d\nu_2$$
whose associated positive measure is 
$$\psi^t\left(
\begin{array}{rl}
\mu_{11} &\mu_{12}\\
\mu_{12} &\mu_{22}
\end{array}
\right)\psi
=\psi_1^2\mu_{11}+2\psi_1\psi_2\mu_{12}+\psi_2^2\mu_{22}.$$

More generally, we will need a notion of  white noise measure with Hilbertian values :

{\bf Definition 2.} {\it Given a bounded positive measure $\mu$ on a measurable space $(E,\mathcal F)$
 and a separable Hilbert space $H$,
we call {\rm $H$-valued white noise measure}  defined thanks to the auxiliary probability space
 $(\Omega,{\mathcal A},\mathbb{P})$ with positive measure $\mu$,   a map from $\mathcal F$ into $L^2((\Omega,{\mathcal A},\mathbb{P});H)$
such that :

($\alpha$) \quad$\forall A\in {\mathcal F},\;\forall h\in H,\quad<\nu(A),h>_H$
is a centered Gaussian variable with variance $\mu(A)\|h\|^2_H$

($\beta$) \quad$A\in{\mathcal F}\longrightarrow \nu(A)$ is $\sigma$-additive in $L^2((\Omega,{\mathcal A},\mathbb{P});H)$

($\gamma$) \quad If $A_1,\ldots,A_k\in{\mathcal F}$ are pairwise disjoint, $\forall h\in H$ the r.v. $<\nu(A_1),h>,\ldots, 
 <\nu(A_k),h>$ are independent.}

Such $\nu$  naturally extends to functions $f\in L^2(E,\mathcal F, \mu)$ and $\forall h\in H$ $<\nu(f),h>_H$ is a centered
Gaussian variable with variance $\int f^2 d\mu\|h\|^2_H$.

To construct such a  $\nu$, let us consider a sequence of independent copies $\nu_n$ of a real white noise measure
  on $(E,\mathcal F)$ with associated positive
 measure $\mu$, and 
for $f\in L^2(E,\mathcal F, \mu)$ let us put
$$\nu(f)=\sum_n \nu_n(f)\chi_n$$
where $\chi_n$ is a complete orthonormal system of $H$.

Similarly to the bivariate case, such a $\nu$ may be transformed in different ways.

a) multiplying by $\varphi\in  L^2(E,\mathcal F, \mu)$
$$(\varphi\nu)(A)=\int1_A\,\varphi\,d\nu$$
$\varphi\nu$ is a $H$-valued white noise measure, s.t.
$\forall f\in L^\infty(E,{\mathcal F}, \mu)$, $\forall h\in H$
$${\mbox{var}}[<(\varphi\nu)(f),h>_H]=\int f^2\varphi^2\,d\mu.\|h\|_H^2.$$
\indent b) For $x\in H$, we can define the scalar white noise measure $(x,\nu)$ by
$$(x,\nu)(A)=<x,\nu(A)>_H.$$
whose associated positive measure is $\|x\|^2\mu$.

c) For $\psi\in L^2((E,{\mathcal F}, \mu);H)$ we can define the scalar white noise measure $(\psi,\nu)$ with associated positive measure
$\|\psi\|_H^2.\mu$ in the following way :

if $\psi$ is decomposable $\psi=\sum_{i=1}^k \psi_i(w).h_i$ then we put
$$(\psi,\nu)=\sum_{i=1}^k\psi_i(w).(h_i,\nu)$$
where $(h_i,\nu)$ is defined in b). The associated positive measure is
$$
\sum_{ij}\psi_i(w)\psi_j(w)<h_i,h_j>_H.\mu=\|\sum_{i=1}^k\psi_i(w)h_i\|_H^2.\mu.$$
For the general case, let $\psi_n$ be decomposable s.t. $\psi_n\rightarrow\psi$ in $L^2((E,{\mathcal F}, \mu);H)$, then we put
$$(\psi,\nu)(A)=\lim_n(\psi_n,\nu)(A){\mbox{ in }}L^2(\Omega,{\mathcal A},\mathbb{P}).$$

After these preliminaries, we can  propose an answer to our initial question: let us consider an error structure 
$S=(W,\;{\mathcal W},\;m,\;\mathbb{D},\;\Gamma)$ admitting a D-gradient $D$ constructed with the help of the separable Hilbert space $H$.

To any $X\in \mathbb{D}$ we shall associate  a real white noise measure that
 will be called its {\bf measure-valued gradient} with satisfactory properties by image. 

\noindent{\bf Definition 3.} {\it  Let $\nu$ be an $H$-valued white noise measure on $(W,\mathcal F)$ with associated positive measure $m$. 
Let $X\in \mathbb{D}$, and let $DX$ be
 its D-gradient constructed with the Hilbert space $H$. The scalar white noise measure
$$(DX,\nu)$$
defined as in c) above, will be called the measure-valued gradient of $X$ and denoted $$d_G X.$$
Thus $\forall f\in L^2(W,{\mathcal W},m)$ we have
$$d_G X(f)\;\left(=\int fd_G X\right)\;=<DX,\nu(f)>_H.$$
 Similarily if $X\in\mathbb{D}^d$ its measure-valued gradient is defined as the column-vector of the measure-valued gradients of its components.
It is therefore an $\mathbb{R}^d$-valued  white noise measure\footnote{The $G$ of $d_G X$ is for Gauss who may be
 considered as the founder of error propagation calculus cf [1].}.} 

\noindent{\bf Proposition 1.}{ \it Let $X\in\mathbb{D}$. Let us denote 
$$\mathbb{E}_{\mathbb{P}}(d_GX)^2$$
the associated positive measure of $d_GX$. We have $\mathbb{E}_{\mathbb{P}}(d_GX)^2<<m$ and 
$$\frac{\mathbb{E}_{\mathbb{P}}(d_GX)^2}{dm}=\Gamma[X].$$}

\noindent{\bf Proof.} Let $f\in L^\infty(W,{\mathcal W}, m)$. The Gaussian r.v. $\int f d_GX=<DX,\nu(f)>_H$ is defined on $(\Omega, {\mathcal A},\mathbb{P})$ 
and has as variance $\int f^2<DX,DX>_Hdm$ by the construction c) defining $(DX,\nu)$. Hence, $\mathbb{E}_{\mathbb{P}}(d_GX)^2$ has a density
with respect to $m$ equal to $\Gamma[X]$. \hfill Q.E.D.

Similarly, if $X\in\mathbb{D}^d$
$$\mathbb{E}_{\mathbb{P}}(d_GX(d_GX)^t)=\underline{\underline{\Gamma}}[X,X^t].m.$$
where $\underline{\underline{\Gamma}}[X,X^t]$ is the matrix with elements $\Gamma[X_i,X_j]$.

\noindent{\bf Proposition 2.} {\it

\noindent a) $\forall X\in\mathbb{D},\quad\forall F\in Lip$
$\quad\quad d_G(F\circ X)=F^\prime(X)d_G X$

\noindent b) $\forall X\in\mathbb{D}^d,\quad\forall F\in{\mathcal C}^1\cap Lip(\mathbb{R}^d)$
$\quad d_G(F\circ X)=\sum_{i=1}^d F^\prime_i(X)\,d_G X_i.$}

\noindent{\bf Proof.} These properties are straightforward from the corresponding ones of the D-gradient.
\section{Images}
Let us look now at what happens by image. Let $X\in\mathbb{D}^d$ and let \break $S_X=(\mathbb{R}^d,\; {\mathcal B}(\mathbb{R}^d),\;
 X_\ast m,\; \mathbb{D}_X,\; \Gamma_X)$ be the image by $X$ of the error structure $S=(W,\;{\mathcal W},\;m,\;\mathbb{D},\;\Gamma)$.

For $F\in \mathbb{D}_X$, to define $d_GF$ we put $$d_GF=X_*(d_G(F\circ X))$$ i.e. $d_GF$ is the image by $X$ of the white noise measure 
$d_G(F\circ X)$\break $=(D(F\circ X),\nu).$
It is defined as a usual image of measure by $$(d_GF)(A)=(d_G(F\circ X))(X^{-1}(A))\quad\forall A\in {\mathcal B}(\mathbb{R}^d)$$
or $\int u^t.d_GF=\int u^t\circ X.d_G F\circ X$ for $u\in L^\infty(X_*m,\mathbb{R}^d)$.

Similarly if $\Phi=(\Phi_1,\ldots,\Phi_k)\in\mathbb{D}^k$, $d_G\Phi$ is defined as the column vector $(d_G\Phi_i)$.

\noindent{\bf Proposition 3.} {\it 

a) $\forall F\in \mathbb{D}_X$ the positive measure associated to $d_GF$ is absolutely continuous w.r. to $X_*m$ and 
$$\frac{\mathbb{E}_{\mathbb{P}}(d_GF)^2}{dX_*m}=\Gamma_X[F].$$ 

b) $\forall F\in {\mathcal C}^1\cap Lip(\mathbb{R}^d)$ $$d_GF=\nabla F^t.d_GI$$ where $\nabla F$ is the usual gradient of $F$ and $I$ is the identity map from $\mathbb{R}^d$
onto $\mathbb{R}^d$which belongs to $(\mathbb{D}_X)^d$. The $\mathbb{R}^d$-valued white noise measure $d_GI=X_*d_GX$ has for associated positive matrix of measures
$(\Gamma_X[I_i,I_j].X_*m)_{ij}$.}

{\bf Proof.} a) From the fact that $\mathbb{E}_{\mathbb{P}}[(d_G(F\circ X))^2]=\Gamma[F\circ X].m$, the image of the white noise measure $d_G(F\circ X)$ by $X$ has for associated
positive measure $\mathbb{E}[\Gamma[F\circ X]|X=x].X_*m$ because of the definition of the conditional expectation, i.e. $\Gamma_X[F].X_*m$. This part of the proposition 
shows that the property of proposition 1 is preserved by image. 

b) We know by proposition 2 that if $ F\in {\mathcal C}^1\cap Lip$ $$d_G(F\circ X)=(\nabla F)^t\circ X\;d_GX$$ the result follows taking the image.\hfill Q.E.D.

Let us denote $L^2(\mathbb{R}^d, \underline{\underline{\Gamma}}_X[I].X_*m)$ the space of $d$-uples of functions\break $v=(v_1,\ldots,v_d)$ defined on
 $(\mathbb{R}^d,{\mathcal B}(\mathbb{R}^d))$ equipped with the norm given by $\| v\|^2=\int v^t\underline{\underline{\Gamma}}_X[I]v\;dX_*m.$ We obtain the main result of 
our study :

\noindent{\bf Proposition 4.} {\it For every $F\in\mathbb{D}_X$, there exists an element of $L^2(\mathbb{R}^d, \underline{\underline{\Gamma}}_X[I].X_*m)$ denoted 
$\nabla_XF$
such that $$d_GF=(\nabla_XF)^td_GI.$$ We have $\underline{\underline{\Gamma}}_X[F]=(\nabla_XF)^t\underline{\underline{\Gamma}}_X[I]\nabla_XF$ and on the
initial structure we have also $$d_G(F\circ X)=(\nabla_XF)^t\circ X\;d_GX.$$}

\noindent{\bf Proof.} Let $(F_n)$ be a sequence of ${\mathcal C}^1\cap Lip$ functions converging to $F$ in $\mathbb{D}_X$. Denoting ${\mathcal E}_X$ the Dirichlet form
of the structure $S_X$, we have
$${\mathcal E}_X[F_n-F_m]=\frac{1}{2}\int\nabla(F_n-F_m)^t\underline{\underline{\Gamma}}_X[I]\nabla(F_n-F_m)\;dX_*m$$ hence $\nabla F_n$ converges 
in $L^2(\mathbb{R}^d, \underline{\underline{\Gamma}}_X[I].X_*m)$. Its limit $\xi$ doesn't depend on the used sequence. For all $u\in L^\infty(\mathbb{R}^d, X_*m)$
the Gaussian variables $\int u\;d_GF_n=\int u(\nabla F_n)^t\;d_GI$ converge in $L^2(\Omega, {\mathcal A},\mathbb{P})$ to $\int u\,\xi^t\,d_GI$ in other words the Gaussian variables
$\int u\circ X\,d_GF_n\circ X$ converge in $L^2(\Omega, {\mathcal A},\mathbb{P})$ to $\int u\circ X\,\xi^t\circ X\;d_GX$. It follows that $d_GF=\xi^t\,d_GI$ and
$d_G(F\circ X)=\xi^t\circ X\,d_GX$.\hfill Q.E.D.

For a function $U\in(\mathbb{D}_X)^p$ with values in $\mathbb{R}^p$ we denote $\nabla_XU$ the matrix of the $\nabla_X$ of its components. We obtain a differential calculus :

\noindent{\bf Proposition 5.} {\it Let $U$ be a map from $\mathbb{R}^d$ into $\mathbb{R}^p$ such that $U\in(\mathbb{D}_X)^p$ and $V$ a map from $\mathbb{R}^p$ into $\mathbb{R}^q$ such that
$V\circ U\in (\mathbb{D}_X)^q$ and $V\in(\mathbb{D}_{U\circ X})^q$. Then
$$(\nabla_X(V\circ U))^t=(\nabla_{U\circ X} V)^t\circ U.(\nabla_XU)^t.$$}
\noindent{\bf Proof.} By proposition 4 applied to $U$ we have
$$d_GU=(\nabla_XU)^t\,d_GI_d,\quad\quad d_G(U\circ X)=(\nabla_X U)^t\circ X\,d_GX,$$
by proposition 4 applied to $V\circ U$ we have
$$d_G(V\circ U)=(\nabla_X(V\circ U))^t\,d_GI_d,\quad\quad d_G(V\circ U\circ X)=(\nabla_X(V\circ U))^t\circ X\,d_GX,$$
now by proposition 4 applied to $V$ on the image structure by $U\circ X$ we have
$$d_GV=(\nabla_{U\circ X})^t\,d_GI_p,\quad\quad d_G(V\circ U\circ X)=(\nabla_{U\circ X}V)^t\circ U\circ X\,d_G(U\circ X).$$ It follows that
$$(\nabla_X(V\circ U))^t\circ X=(\nabla_{U\circ X}V)^t\circ U\circ X.(\nabla_X U)^t\circ X$$ equality in the space $L^2(E,{\mathcal F}, 
\underline{\underline{\Gamma}}[X].m)$ and 
$$(\nabla_X(V\circ U))^t=(\nabla_{U\circ X}V)^t\circ U.(\nabla_X U)^t$$  equality in the space $L^2(\mathbb{R}^d, \underline{\underline{\Gamma}}_X[I].X_*m)$. The 
argument comes therefore from the fact that the notions are defined thanks to images of measures.\hfill Q.E.D.

Let $M_X$ be a measurable square root (non necessarily positive) of the matrix  $\underline{\underline{\Gamma}}_X[I]$, i.e. such that 
$M_X^tM_X=\underline{\underline{\Gamma}}_X[I]$.

\noindent{\bf Corollary} {\it For $F\in \mathbb{D}_X$ let us define $$D_X F=(\nabla_XF)^t M_X^t$$
then $D_X$ is a Dirichlet-gradient for the image structure $S_X$ defined with the Hilbert space $\mathbb{R}^d$.}

\noindent{\bf Proof.} $(D_XF,D_XF)_{\mathbb{R}^d}=(\nabla_XF)^t\underline{\underline{\Gamma}}_X[I]\nabla_XF$ which is equal to
$\underline{\underline{\Gamma}}_X[F]$ by proposition 4. Hence $D_X$ is a D-gradient for $S_X$.\hfill Q.E.D.

\section{Example}
Let us consider the classical Wiener space $(W,{\mathcal W},m)$ with $W={\mathcal C}_0[0,1]$, ${\mathcal W}$ its Borel $\sigma$-field and $m$ the Wiener measure
equipped by the Ornstein-Uhlenbeck structure $(W,{\mathcal W},m, \mathbb{D},\Gamma)$ characterized by
$$\forall f\in L^2[0,1],\qquad\Gamma[\int f\,dw]=\|f\|^2_{L^2}$$ the space $\mathbb{D}$ is usually denoted $D_{2,1}$ or $\mathbb{D}_1^2$ (cf [4], [7], [8], [9]). We consider the
Feyel-la-Pradelle gradient denoted $\#$, it is a linear map from $\mathbb{D}$ into $L^2(m,L^2(\hat{W},\hat{\mathcal W}, \hat{m}))$ where $(\hat{W},\hat{\mathcal W}, \hat{m})$
is a copy of $(W,{\mathcal W},m)$. Thanks to the functional calculus
it is characterized by its values on the first chaos :
$$\forall f\in L^2[0,1],\qquad(\int f\,dw)^\#=\int f\,d\hat{w}$$
the Hilbert space $H$ is $L^2(\hat{W},\hat{\mathcal W}, \hat{m})$. Let $(Z_n)$ be an orthonormal basis of $L^2(W,{\mathcal W},m)$ for instance composed with
 a basis of each Wiener chaos, $(\hat{Z_k})$ the corresponding basis of $L^2(\hat{W},\hat{\mathcal W}, \hat{m})$ and let $g_{n,k}$ be i.i.d. reduced 
Gaussian variables defined on a probability space $(\Omega, {\mathcal A}, \mathbb{P})$. 

Putting for $Y\in L^2(W,{\mathcal W},m)$ 
$$\int Y\,d\nu= \sum_{k,n}\mathbb{E}_m[YZ_n]\hat{Z_k}\,g_{n,k}$$ defines according to definition 2 an H-valued white noise measure on $(W,{\mathcal W})$ with positive
measure $m$.

If $X\in \mathbb{D}$ according to definition 3 for $Y\in L^\infty(W,{\mathcal W},m)$
$$\int Y d_GX=\sum_{k,n}\mathbb{E}_m[YZ_n\mathbb{E}_{\hat{m}}[X^\#\hat{Z_k}]]g_{n,k}$$ and we have
$$\begin{array}{rl}
\mathbb{E}_{\mathbb{P}}[(\int Y\,d_GX)^2]&=(\sum_{k,n}\mathbb{E}_m[YZ_n\mathbb{E}_{\hat{m}}[X^\#\hat{Z_k}]])^2\\
&=\sum_k\mathbb{E}_m[(Y\mathbb{E}_{\hat{m}}[X^\#\hat{Z_k}])^2]=\mathbb{E}_m[Y^2\Gamma[X]]
\end{array}
$$
so that the positive measure associated with $d_GX$ is indeed $\Gamma[X].m$ and the study applies.

These results mean that a differential calculus may be defined on an error structure and its images, 
satisfying the expected coherence property, which coincides with the usual differential calculus on ${\mathcal C}^1\cap Lip$
functions but exists also by completion for any function in the Dirichlet spaces of the images structures, coherence
 being preserved, thanks to the fact that the image of a gradient is now defined as the usual image of a measure.  The tools introduced here are 
not intrinsic, this would be an interesting program to geometrize them. But
in the applications, for studying the sensitivity of stochastic models, we are mostly concerned with computations in  situations where an error structure is 
defined on the Wiener space (e.g. the Ornstein-Uhlenbeck structure or a generalized Mehler-type structure) or on the Poisson space or
both, and all is about images of this structure (cf. [3]) the preceding study is relevant from this point of view.



\end{document}